\overfullrule=0pt
%
%

%
%

%


\input mssymb
%
%
\font\euler=eufm10
\font\eulers=eufm7
\font\eulerss=eufm5

\newfam\eulerfam
\textfont\eulerfam\euler
\scriptfont\eulerfam\eulers
\scriptscriptfont\eulerfam\eulerss

\def\frak#1{{\fam\eulerfam\euler #1}}

\def\eqdef{\buildrel \rm def \over =}
\def\rest{\mathord{\restriction}}

\def\bb{\bigbreak}

\def\phi{\varphi}
\def\cov{{\rm cov}}

\def\sat{\models}
\def\su{\subseteq}
\def\a{\alpha}
\def\b{\beta}

\def\e{\nu}
\def\d{\delta}
\def\l{\lambda}
\def\k{\kappa}

\def\om{\omega}

\def\lng{\langle}
\def\rng{\rangle}
\def\ov{\overline}
\def\sm{\setminus}
\def\nac{{\rm nacc }\,}
\def\nacc{\nac}
\def\acc{{\rm acc}\,}

\def\cf{{\rm cf}}
\def\inv{{\rm Inv}}

\def\epsilon{\nu}
\def\otp{{\rm otp}\,}
\baselineskip=16pt
\def\ran{{\rm  ran}}

\def\endproof#1{\hfill $\smiley_{#1}$}

\def\imply{\Rightarrow}
\def\iff{\Leftrightarrow}

\def\id{{\rm id}}
\def\proof{\smallbreak\noindent{\sl Proof}: }
\def\Inv{{\rm INV}}

\def\tp{{\rm tp}}

\def\conc{\concatenate}
\def\concatenate{\widehat{\;}}
\def\univ{{\rm Univ}}

\def\ba{{\bf a}}
\def\bb{{\bf b}}
\def\bc{{\bf c}}
\def\init{\triangleleft}
\def\forks{\fork}
\baselineskip=27pt
\def\bbM{{\frak C}} 
\def\bI{{\bf I}}
\def\bJ{{\bf J}}
\def\av{{\rm Av}}
\def\aut{{\rm Aut}}


\newbox\noforkbox \newdimen\forklinewidth
\forklinewidth=0.3pt   
\setbox0\hbox{$\textstyle\bigcup$}
\setbox1\hbox to \wd0{\hfil\vrule width \forklinewidth depth \dp0
			height \ht0 \hfil}
\wd1=0 cm
\setbox\noforkbox\hbox{\box1\box0\relax}
\def\unionstick{\mathop{\copy\noforkbox}\limits}
\def\nonfork#1#2_#3{#1\unionstick_{\textstyle #3}#2}
\def\nonforkin#1#2_#3^#4{#1\unionstick_{\textstyle #3}^{\textstyle #4}#2}     

\setbox0\hbox{$\textstyle\bigcup$}
\setbox1\hbox to \wd0{\hfil$\nmid$\hfil}
\setbox2\hbox to \wd0{\hfil\vrule height \ht0 depth \dp0 width
				\forklinewidth\hfil}
\wd1=0cm
\wd2=0cm
\newbox\doesforkbox
\setbox\doesforkbox\hbox{\box1\box0\relax}
\def\nunionstick{\mathop{\copy\doesforkbox}\limits}

\def\fork#1#2_#3{#1\nunionstick_{\textstyle #3}#2}
\def\forkin#1#2_#3^#4{#1\nunionstick_{\textstyle #3}^{\textstyle #4}#2}

\font\circle=lcircle10

\setbox0=\hbox{~~~~~}
\setbox1=\hbox to \wd0{\hfill$\scriptstyle\smile$\hfill} 
\setbox2=\hbox to \wd0{\hfill$\cdot\,\cdot$\hfill} 
\setbox3=\hbox to \wd0{\hfill\hskip4.8pt\circle i\hskip-4.8pt\hfill} 

\wd1=0cm
\wd2=0cm
\wd3=0cm
\wd4=0cm

\newbox\smilebox
\setbox\smilebox \hbox {\lower 0.4ex\box1
		 \raise 0.3ex\box2
		 \raise 0.5ex\box3
		\box4
		\box0{}}
\def\smiley{\leavevmode\copy\smilebox}


\newcount\secno
\newcount\theono   

\catcode`@=11
\newwrite\mgfile

\openin\mgfile \jobname.mg
\ifeof\mgfile \message{No file \jobname.mg}
	\else\closein\mgfile\relax\input \jobname.mg\fi
\relax
\openout\mgfile=\jobname.mg

\newif\ifproofmode
\proofmodefalse            

\def\@nofirst#1{}

\def\neusection{\advance\secno by 1\relax \theono=0\relax}
\def\neuchap{\secno=0\relax\theono=0\relax}

\neuchap

\def\labelit#1{\global\advance\theono by 1%
             \global\edef#1{%
             \number\secno.\number\theono}%
             \write\mgfile{\@definition{#1}}%
}


\def\ppro#1#2:{%
\labelit{#1}%
\smallbreak\noindent%
\@markit{#1}%
{\bf \ignorespaces #2:}}

\def\@definition#1{\string\def\string#1{#1}
\expandafter\@nofirst\string\%
(\the\pageno)}

\def\@markit#1{
\ifproofmode\llap{{\tt \expandafter\@nofirst\string#1\ }}\fi%
{\bf #1\ }
}
 
\def\labelcomment#1{\write\mgfile{\expandafter
		\@nofirst\string\%---#1}} 

\catcode`@=12


\centerline{\bf The Universality Spectrum of Stable Unsuperstable
Theories}

\centerline{By Menachem Kojman and Saharon
Shelah\footnote\dag{Partially supported by the Binational Science
Foundation. Publication No. 447}}

\centerline{Mathematics Department, The Hebrew University of
Jerusalem, Jerusalem} 

\centerline {ABSTRACT}

It is shown that if $T$ is stable unsuperstable, and
$\aleph_1<\l=\cf\l<2^{\aleph_0}$, or
$2^{\aleph_0}<\mu^+<\l=\cf\l<\mu^{\aleph_0}$ then $T$ has no universal
model in cardinality $\l$, and if e.g. $\aleph_\om<2^{
\aleph_0}$ then
$T$ has no universal model in $\aleph_\om$. These results are
generalized to $\k=\cf\k<\k(T)$ in the place of $\aleph_o$. Also: if
there is a universal model in $\l>|T|$, $T$ stable and $\k<\k(T)$ then
there is a universal tree of height $\k+1$ in cardinality $\l$.

\neusection
{\bf \S\number\secno. Introduction}
 
We handle the universal spectrum of stable-unsuperstable first order
theories. This continues [KjSh 409] and adds information the picture
started
up in [Sh 100]. The general subject addressed here is the universal
model problem, which although natural and old, was not treated very
extensively in the past. For
background, motivation and history of the subject see the introduction
to [KjSh~409], a paper in which unstable theories with the strict
order property are handled (e.g. the class of linear orders).

When looking at a class $K$  of structures together with a class of allowed
embeddings --- say all models of some 
first order 
theory (T) with elementary embeddings --- we get a partial order:
$A\le B$ if there is a mapping of $A$ into $B$ in the class of allowed
mappings. 
The universal model problem can be phrased, in this context, as a
question about this partial order: is there in $\{M\in K:||M||\le\l\}$
 a ``greatest''
element --- which we call ``universal'' --- namely one  such that
all other elements $M\in K$, $||M||\le\l$   are smaller than or equal
to it. This question can 
be elaborated:
what is the cofinality, i.e.
 the minimal cardinality of  a subcollection of
elements such that 
every element is smaller than or equal to {\sl one} of the elements in
this 	
subcollection?
 Can a universal object be found outside our
collection?
 (for instance, is there a model of $T$ of cardinality
$\mu>\l$ such that every model of $T$ of cardinality $\l$ is
elementarily embeddable into it). How does the existence, or
nonexistence, of universal objects in one collection of structure
influence the existence or non existence of universal objects in
related collections?

 In this paper we prove that if $T$ is stable unsuperstable, and
$\mu^+<\l=\cf\l<\mu^{\aleph_0}$ then $T$ has no  model
of cardinality $\l$ into which all models of $T$ of the same
cardinality are elementarily embeddable, not even a family of models $\lng
M_i:i<\l'\rng$, $\l'<\mu^{\aleph_0}$,  each of cardinality
$\l$ such that every model of $T$ of 
cardinality $\l$ is elementarily embedded into some model in the
family. It follows from the theory of covering numbers
that certain 
singular cardinals are also not  in the universality spectrum of
stable unsuperstable theories.

 Also, it is shown that a certain theory (the ``canonical'' stable
unsuperstable theory)
is ``minimal'' with respect to the existence of universal models,
namely that whenever some stable unsuperstable theory $T$ has a
universal model in cardinality $\l$, also this theory has one.

 	We
mention here without proof 
that $GCH$ implies 
that all first order theories have universal models in all uncountable
cardinals (above the cardinality of the theory), and that the question
whether $\aleph_1$ is in the universality spectrum of a countable,
stable but not superstable theory is independent of $ZFC+ 2^{\aleph_0}
=\aleph_2$ (see [Sh 100] \S2). At this point it is interesting to note
that it is consistent that there is a universal graph in
$\aleph_2<2^{\aleph_0}$ (see [Sh 175a]),
  but it is not consistent to have
a universal model for some countable, stable unsuperstable
$T$. So in this respect, stable unsuperstable theories
are not ``$\le$'' all unstable theories. 

In  subsequent papers, 
 universality spectrums of   some classes 
of infinite abelian groups, and complementary consistency results to
the negative results known so far will be dealt with. (Note: if $T$ is
countable first order, stable unsuperstable theory, and
$\mu=\Sigma_n\mu_n^{\aleph_0}$, then there is a universal model for
$T$ in $\mu$; if, say, $\beth_\om^+<\mu<\beth_{\om+1}$ there isn't; and
we do not settle here the case $\mu=\beth_\om^+$.)

We assume some familiarity  with the definitions
of stability and superstability, as well as with  fundamentals of
forking theory (to be found  in e.g. [ [Sh-c],III).

Lastly, those reader who speculate that \smiley\footnote{*}{We thank
Martin Goldstern for the smiley \TeX nology and for
$\nonfork{}{}_{}$.} is some inexplicable 
whim of the laser printer are wrong. Smileys indicate ends of proofs,
and replace the old and square boxes. 
  
\neusection 

{\bf \S\number\secno. Preliminaries and Setup}

Having fixed attention on a given $T$, we work in some ``monster
model'' $\bbM$, which is a big saturated model, of which all the
models we are interested in are elementary submodels of smaller cardinality.

\ppro\firstdefs Definition: for a complete first order theory $T$,
\item{(0)} A model $M\sat T$ is {\sl $\prec$-universal in cardinality} $\l$ if
$||M||=\l$ and for every $N\sat T$ such that $||N||=\l$ there is an
elementary embedding $h:N\to M$. It is $<$-universal, if we omit
``elementary'' from the definition. 
\item {(1)} $\univ(T,\prec)=\{\l:\l\;\;{\rm is}\;\;{\rm a}\;\;{\rm
cardinal}\;\;{\rm and}\;\; T\;\; {\rm has\;\; a\;\; universal \;\;
model \;\; in}\;\;\l\}$ is the {\it universality spectrum} of $T$.

\item{(2)}  $\univ_p (T, \prec)$  is the family of pairs $(\l,\mu)$
such that 
there is a family of $\mu$	  models of $T$
each of cardianlity $\l$, such that any model of $T$ 
of cardianlity $\l$ can be elelmentarily embeded
into one of them.
\item{(3)}  $\univ_t (T, \prec )$  is the family of triples 
$(\l,\k, \mu )$ such that there is is a family of $\mu$
models of $T$ each of cardianlity $\le  \k$ and any model of $T$
of cardianlity $\l$ can be elemntarily embeded into one
of them.

In this paper ``universal'' means $\prec$-universal and $\univ(T)$
means $\univ(T,\prec)$ unless otherwise stated.

\ppro\seconddefs Definition:                                              

\item{(0)} A Theory $T$ is {\it stable in $\l$} if for every model $N$ and
set  $A\su N$, $|A|\le \l\imply |S(A)|\le \l$. For equivalent
definitions see [Sh-g ],II.2.13.
\item {(1)} $\k(T)$, the 
cardinal of $T$ is as defined in [Sh-c],III,\S3.
We
recall from [Sh-c]III that for a countable, complete first
order $T$, $T$ is stable unsuperstable iff $\k(T)=\aleph_1$.
\item {(2)} The notation $\nonforkin \ba A _ B ^N$ means ``the type of
$\ba$ over the set $A$ in the model $N$ does not fork over the set
$B$''. The notation $\forkin \ba A _ B ^N$ means  ``the type of
$\ba$ over the set $A$ in the model $N$ forks over the set $B$''. When
the model $N$ in which the relation of forking exists is clear from
context, it is omitted.

By small bold faced letters we shall denote finite sequences of
elements from a model. Following a widely spread abuse of notation we
shall not write $\ba\in |N|^{<\om}$, but write $\ba\in N$, and even
refer to $\ba$ as ``element''. This is perfectly all right with what is
 about to be done here, because we may add the finite sequences as
elements into
the model and work in $T^{eq}$ (or $\bbM^{eq}$), or replace a type of an $n$-tuple by a
1-type when necessary.

The forking facts which we shall  need  are summarized in the
following quotation from [Sh-c] p.84:
\bigbreak
\ppro \forkFacts Theorem:
\item{(0)} (finite character of forking) If $\fork\ba A _ B$ then there
is some {\it finite} set $A'\su A$ such that $\forks \ba {A'} _ B$.
Also   $\fork \ba A _ B$  iff  $\fork  \ba  {A \cup B}  _ B$.
\item{(1)} (symmetry) $\nonfork {\ba}{A\cup\bb}_A$ iff $\nonfork
{\bb}{A\cup \ba}_A$.
\item {(2)} (transitivity) if $A\su B\su C$ and $\nonfork {\ba}{C}_B$
and $\nonfork \ba B _ A$ then $\nonfork \ba C _ A$
\item {(3)} Let $B\su A$ , then: $\nonfork \bb {A\cup\ba} _ {B\cup A}$
and $\nonfork \ba A _ B$ iff $\nonfork {\ba\conc\, \bb} A _ B$
\item{(4)} When $M$ is a model, the type $p$ does not fork over $|M|$
iff $p$ is finitely satisfiable in $M$.
\item {(5)} if $A\su B\su C$, $p\in S(B)$ does not fork over $A$, then
there is some $q\in S(C)$,$p\su q$ and $q$ does not fork over $A$. 
\item {(6)}([Sh-c] p.113) If $p\in S(|M|)$ is definable over $A$ where
$A\su M$ then $p$ does not fork over $A$.
 
We need a few facts about sets of indiscernibles. We denote sets of
indiscernibles by $\bI$ and $\bJ$. We say that
$\tp(\bI)=\tp(\bJ)$ if for every $n$, formula $\phi$ and elements
$\ba_1,\cdots,\ba_n\in \bI,\bb_1,\cdots,\bb_n\in \bJ$,
$\tp(\ba_1,\cdots,\ba_n,\emptyset)=\tp(\bb_1,\cdots,\bb_n,\emptyset)$. 

\ppro \indisFacts Theorem:
\item {(1)} ([Sh-c],III,4.13, p.77) If $T$ is stablle, $\phi(\ov  x,\ov y)$ a
formula, then there is some natural number $n(\phi)$ such that for
every set of indiscernibles $\bI$ and parameters $\bc$,  either
$|\{\ba\in \bI:\sat \phi(\ba,\bc)\}|<n(\phi)$ or
$|\{\ba\in\bI:\sat\neg\phi(\ba,\bc)\}|<n(\phi)$ 
\item {(2)} ([Sh-c],III, 1.5 p.89) Let $I$ be an infinite set of indicernibles.
$\av_\Delta(\bI,A)$, the average of 
$\bI$ over the set of formulas $\Delta$ and over the set of parameters $A$, is the
set of all fomulas $\phi(\ov x,\bc)$ such that $\phi\in \Delta$,
$\bc\in A$ and $\sat \phi(\ba,\bc)$ for all but finitely many $\ba\in
\bI$.
\item{(3)}([Sh-c],III,, 3.5 p. 104) If $\bJ$ is an indiscernible set
over $A$, $B$ is any set, then there is $\bI\su\bJ$ such that $\bJ -
\bI$ is indiscernible over $A\cup B\cup \bigcup \bI$, and
\item\item{(a)} $|\bI|
\le \k(T)+|B|$.
\item\item{(b)} If $|B|<\cf(\k(T))$ then $
|\bI|<\k(T)$. (The interesting
case is when $|\bJ|$ is large enough with relation to $|B|$.)
\item {(4)}([Sh-c],III,4.17  p.117) If $\bI,\bJ$ are infinite indiscernible sets,
and $\av(\bI,\bigcup\bI)=\av(\bJ, \bigcup\bI)$ and
$\av(\bJ,\bigcup\bJ)=\av(\bI,\bigcup\bJ)$ then $\av(\bI,\bbM)=\av(\bJ,\bbM)$.
\item{(5)} ([Sh-c] p.112),III, 4.9 If $\Delta$ is finite and $p\in S^m(|M|)$,
then for every type $q\in S^m(B)$ extending $p$ which does not fork
over $M$ there is an
infinite $\Delta$-indiscernible set $\bI\su M$ such that $q=\av_\Delta(\bI,B)$.
\item{(5)}([Sh-c],III,1.12 p.92) For every  $\b$ and set $A$ there is
an indiscernible sequence $\bI$ over $A$ and based on $A$ (=  for
every $B$ $\av(\bI,B)$ does not fork over $A$) such that $\bb\in
\bI$.
\vskip1cm

The interested reader is welcome to inquire [Sh-c] for more details
and/or results.

We recall some combinatorics which we need:
\ppro\guessDef Definition:

 Suppose $\l$ is a regular uncountable cardinal, and $S\su \l$
is stationary.
\item{(1)} A sequence $\ov C=\lng c_\d:\d\in S\rng$ is a {\it club
guessing sequence on $S$} if $c_\d$ is a club (=closed unbounded
subset) of $\d$ for every $\d\in S$
and for every club $E$ of $\l$ the set $S_E=\{\d\in S:c_\d\su E\}$ is
stationary.
\item {(2)} For $\ov C$ as in (1), $\id^a(\ov C)\eqdef\{A\su S
:$
there is a club $E\su\l$ such that $\d
\in
 A\cap S\imply c_\d\not\su E\}$ is a $\l$-complete proper ideal.
\item {(3)} a sequence $\lng P_\d:\d\in S\rng$, $S\su \l$, is a {\it
weak } club 
gueesing sequence if $P_\d=\lng c^\d_i:i<i(\d)\rng$, $i(\d)\le \l$,
for each $i<i(\d)$, 
$c^\d_i$ is a club of $\d$ and for every club $E\su
\l$, the set $S_E=\{\d\in S:(\exists i<i(\d))(c^\d_i\su E)\}$ is
stationary. The existence of a weak club guessing sequence is clearly
equivalent to the existence of a sequence $\lng c_\b:\b<\l\rng$ such that
$c_\d\su \b$ and for every club $E\su \l$ the set $\{\a<\l:(\exists
\b)(\sup c_\b=\a)\;\&\; c_\b\su E\}$ is stationary. We call such  a
sequence also a weak club guessing sequence.
\item {(4)} If $\ov P = \lng P_\d:\d\in S\rng$ is a weak club guessing
sequence, then $\id^a(\ov P)=\{A\su S: (\exists E)(E\su \l$ 
is club 
such that  $(\forall \d\in E\cap S)(\neg\exists
i<i(\d))(c^\d_i\su E)\}$ is a proper $\l$-complete ideal. 
\endproof{\guessDef} 

\ppro\guessfacts Fact: 
\item {(1)} If $\l=\cf\l>\aleph_1$ then there are a stationary $S\su
\l$ and
a club guessing sequence  $\ov
C=\lng c_\d:\d\in S\rng$ on $S$ such that for every $\d\in S$ the order
type of $c_\d$ is $\om$.
\item{(2)} If $\k$ is regular and uncountable, $\k^+<\l=\cf\l$, then
there are sequences $\ov C=\lng c_\d:\d\in S\rng$, $S\su \l$
stationary, and $\lng P_\a:\a\in \l\rng$ such that $\otp c_\d=\k$,
$\sup c_\d=\d$, $\ov C$ is a club guessing sequence, $|P_\a|<\l$ and
for every $\d\in S$ and $\a\in \nacc c_\d$, $c_\d\cap \a\in P_\a$. 
\item{(3)} Suppose  $\mu^+<\l=\cf\l$ and $\cf\mu\le\mu$. Then there
is  a weak club guessing sequence  $\ov
C=\lng c_\b:\b<\l\rng$ such that for every $\b<\l$ the order
type of $c_\b$ is $\mu$ and for every $\a< \l$ the set $\{c_\b\cap
\a\}$ has cardinality smaller than $\l$. 

\proof
See [Sh-g], [Sh 365] or the appendix to [KjSh 409] for a proof of (1)
and see [Sh 420]\S1 for the proofs of (2) and (3).\endproof{\guessfacts}

On covering numbers see [Sg-g]. We refer the reader to [KjSh 409],\S4,
for a detailed exposition of 
covering numbers of singular cardinals,
in 
particular to Theorem 4.5 there. Here we quote

\ppro \covDef Definition: $\cov(\l,\mu,\theta,\sigma)$ is the minimal
size of a
family $A\su[\l]^{<\mu}$ which satisfies that for all
$X\in[\l]^{<\theta}$ there are less than $\sigma$ members of $A$ whose
union covers $X$.

\ppro \covFact Theorem:  If $\mu$ is not a fix point of the second
order, i.e. $|\{\l<\mu:\l=\aleph_\l\}|=\sigma<\mu$, 
and $\sigma+\cf\mu<\k<\mu$, then
$\cov(\mu,\k^+,\k^+,\k)=\mu$.\endproof{\covFact}

For example, for every $\aleph_n$ it is true that
$\cov(\aleph_\om,\aleph_{n+1},\aleph_{n+1},\aleph_n)=\aleph_\om$. 
\bigbreak

\neusection

{\bf \S\number\secno. The Machinery}

In this section  $T$ denotes a first order, countable, stable but
unsuperstable theory. 

\ppro \invDef  Definition:
Suppose that $N\sat
T$, and that $\ov N=\lng
N_i:i<\l\rng$ is given, $N_i\prec N_{i+1}\prec N$, $||N_i||<\l$ and
$N_j=\cup_{i<j}N_j$ for limit $j$.  $\ov N$ is called a {\it
representation} of $N$. 

Suppose $c\su \l$ is of limit order type and is enumerated
(continuously and increasingly) by $\lng \a_i:i<\otp c\rng$.
 \item {(0)} for an element $\ba\in N$, $\inv_{\ov
N}(\ba,c)\eqdef\{\a_{i} : \forkin{\ba} {N_{\a_{i+1}}} _ {N_{\a_i}} ^
{N}\}$;
\item{(1)}  $\inv^*_{\ov N}(\ba,c)\eqdef\{{i} : \forkin{\ba}
{N_{\a_{i+1}}} _ {N_{\a_i}} ^ 
{N}\}$
\item{(2)} $P(\ov N,c)\eqdef\{\inv_{\ov N}(\ba,c):\ba\in N\}$
\item{(3)}
$P^*(\ov N,c)\eqdef\{\inv^*_{\ov N}(\ba,c):\ba\in N\}$ 

\ppro \moreInv Definition: Suppose that $\ov C=\lng c_\d:\d \in S\rng$
is a club guessing 
sequence on some staitonary $S\su \l$ and that that $\ov N$ is as in
\invDef.  

\item{(4)}   $\Inv^a (\ov N,\ov C )$  is 
the sequence     $\lng P(\ov N , c_\d) : \d  \in  S\rng$  
modulo the ideal    $\id^a (\ov C)$. 
\item{(4)} Assuming that for all $\d\in S$, $\otp c_\d$ is some fixed
$\d(*)$,  $\Inv^b (\ov N,\ov C)\eqdef\{Y\su \d(*):\{\d\in S:Y\in
P^*(\ov N,c_\d)\}\notin \id^a(\ov C)\}$
\item{(5)} Under the assumptions of (4),        $\Inv^c _{(\ov N,\ov
C)}\eqdef\{Y\su \d(*):\{\d\in S: Y\notin 
P^*(\ov N,c_\d)\}\in \id^a(\ov C)\}$ 

\ppro \moreInvRemark Remark: We shall not use \moreInv\ much,
but our results can be interpreted as saying
that those invariants  do not depend on 
the representation $\ov N$ but  just on the model $N$,
and that  we can prove non universality by just looking at one of
these invariants.

\ppro\InvIsPreservedUnderEmbeddings Lemma:
Suppose $\l=\cf\l>\aleph_1$, $N,M$
are models of $T$, $||N||=||M||=\l$ with given representations $\ov
N,\ov M$. If $h:N\to M$ is an
elementary embedding, then there is some club $E\su \l$ such that for
every
$\ba\in N$ and $c\su E$, $\inv_{\ov N}(\ba,c)=\inv_{\ov M}(h(\ba),c)$.

\proof
 Let $E_h=\{i<\l:\ran(h\rest
N_i)\su M_i\}$. Clearly $E_h$ is a club
of $\l$. Dnote by $N_i^*$ the set $\ran(h\rest N_i)$. So for $\d\in
E_h$, $N_i^*$ is the universe of an elementary submodel of $M_i$.
Denote by $N^*$ the image of $N$ under $h$. 

\ppro\clubClaim Claim:
The set $E_1=\{\d\in E_h: (\forall \ba\in
M_\d)(\nonfork\ba {N^*}_{N_\d(*)})\}$ is a
club.

\proof
As $T$ is countable, stable but unsuperstable, $\k(T)=\aleph_1$. Therefore for
every $\ba\in M$, there is a countable set $A_\ba\su N^*$ such that
$\nonfork\ba{N^*}_{ A_\ba}$. Let $i(\ba)$ be the least 
$  i $
such
that $A_\ba\su N_i^*$. For $\a\in E_h$ let $j(\a)$ be the least 
$j\in E_h$ such
that
for
all $\ba\in [\a]^{<\om}$, 
$i(\ba)\le j$.
$E'=\{\d\in E_h:\a<\d\imply j(\a)<\l\}$ is club. If $\d\in
E'$ and $\ba\in [\d]^{<\om}$, 
then $A_\ba\su N_\d(*)$. So as
$\nonfork\ba{N^*}_{ A_\ba}$, also $\nonfork\ba{N^*}_
{N_\d(*)}$. So $E' \su E_1$. $E_1$ is closed, for if
$\d\in\acc E_1$ and $\ba\in M_\d$, then there is some $\a<\d$
such that $\sup\ba<\a$ and $\a<i<\d$ such that $i\in E_1$.
$\nonfork \ba {N^*}_{N_i^*}$, therefore
$\nonfork \ba {N^*}_{ N_\d(*)}$.\endproof{\clubClaim}

 Let $\lng \a_i:i<\l\rng$ be its increasing
enumeration of $E_1$. We show that for every $\ba\in N$ and $i<\l$,
$\forkin\ba {N_{\a_{i+1}}}_{ N_{\a_i}}^{N} \iff
\forkin {h(\ba)} {M_{\a_{i+1}}}_{ M_{\a_i}}^{M}$. 

 As $E\su E_1$, for every $\a_i$ and $\bb\in M_\d$,
$$\nonforkin{\bb}{N^*}_{N_{\a_i}^*}^{M}$$
 This can be written as
$$\nonforkin{M_{\a_i}}{N^*}_{N_{\a_i}^*}^{M}\leqno{\hskip2cm\star
_1}$$
By monotonicity, for a given $\ba\in N$,
$$\nonforkin{M_{\a_i}}{h(\ba)}_{N_{\a_i}^*}^{M}\leqno\hskip2cm\star_2
$$
Symmetry of non-forking gives
$$\nonforkin{h(\ba)}{M_{\a_i}}_{N_{\a_i}^*}^{M}\leqno\hskip2cm\star_3
$$

Suppose now, first, that
$$\nonforkin{\ba}{N_{\a_{i+1}}}_{N_{\a_i}}^{N}.\leqno$$  As $h$ is
an elementary embedding,
$$\nonforkin{h(\ba)}{N_{\a_{i+1}}^*}_{N_{\a_i}^*}^{M}\leqno\hskip2cm
\star_4$$

By $\star_3$, $\nonfork{h(\ba)}{M_{\a_{i+1}}}_{N_{\a_{i+1}}^*}$ (we
omit
$M$, in which we work from now on). By $\star_4$, and the
transitivity of non forking, 
$\nonfork{h(\ba)}{N_{\a_{i+1}}}_{N_{\a_i}^*}$.
By monotonicity, $\nonfork{h(\ba)}{M_{\a_{i+1}}}_{M_{\a_i}}$.

For the other direction, suppose that
$\nonfork{h(\ba)}{M_{\a_{i+1}}}_{M_{\a_i}}$. By monotonicity,
$\nonfork{h(\ba)}{N_{\a_{i+1}}^*}_{M_{\a_i}}$. By $\star_3$ and the
transitivity of non forking,
$\nonfork{h(\ba)}{N_{\a_{i+1}^*}}_{N_{\a_i}^*}$, which is what we
want.\endproof{\InvIsPreservedUnderEmbeddings}

\ppro \embedCor Corollary: Suppose $\ov N$ and $\ov M$ are as above
and that $h:N\to M$ is an elementary embedding. Let $E$ be the
club given by the previous lemma. If $c\su E$ then
\item{(1)}  for every $\ba\in N$, $\inv_{\ov
N}(\ba,c)=\inv_{\ov M}(h(\ba),c)$;
\item {(2)}  $P(\ov N,c)\su P(\ov M,c)$.

We will need a slight generalization of
\InvIsPreservedUnderEmbeddings:

\ppro \secondEmbedding Lemma:
Suppose $N\sat T$ is with universe $\l$,   $\ov N$ is a representation
of $N$.
Suppose $L\prec M$ are models of $T$, $L$ is of cardinality $\l$,
its universe is $B$ and $\ov L$ is a representation.  If
$h:N\to M$ is an elementary 
embedding, 
then there is some club $E\su \l$ such that for every $c\su E$ and
$\ba\in h^{-1}(B)$,
$\inv_{\ov N}(\ba,c)=\inv_{\ov L} (h(\ba),c)$.

\proof 
Denote by $N_i^*,N^*$ the images of $N_i,N$ under $h$ respectively.
Let $A_i=|N_i^*|\cap B$. Let $A=\cup A_i$. We prove

\ppro \secondClub Claim: There is a club $E_1\su \l$ such that $i\in
E$ implies $\nonfork {N_i^*} A _{A_i}$ and $\nonfork {L_i}{A}_{A_i}$

\proof Same as in \clubClaim.\endproof{\secondClub}

For the rest of the proof, show, precisely as in \clubClaim, that 

$$\nonfork{h(\ba)}{B_{\a_{i+1}}}_{B_{\a_i}}\iff\nonfork
{h(\ba)}{A_{\a_{i+1}}}_{A_{\a_i}}\iff
\nonfork{h(\ba)}{N_{\a_{i+1}}^*}_ {N_{\a_i}^*}$$

When $\lng\a_i:i<\otp c\rng$ is the enumeration of $c$.\endproof{\secondEmbedding}

\ppro \construction Lemma: (the construction Lemma) 
Let $\l$ be
uncountable and regular. Suppose that
 $\ov C$ is a club guessing sequence on some stationary
$S\su \l$ and for every $\d \in S $, $\otp c_\d=\mu$ for some fixed
$\mu$ with $\cf\mu=\aleph_0$. 
  Suppose $Y\su\mu$ is a given set of order type $\om$. Then
there is a model $M\sat T$ of 
cardinality $\l$ and a representation $\ov M$ such that 
for every $\d\in S$,  $Y\in P^*(\ov N,c_\d)$.

\proof
We work in the  monster model, ${\bbM}$. By 
$\k(T)>\aleph_0$, there is some $\bb$ and $M$ with the property that for every
finite set $A$, $\forkin \bb {M}_A ^ {\bbM}$. Pick by induction on $n$ a
finite
sequence $\ba_n$ such that
\item{(i)} $\ba_n$ is a proper initial segment of $\ba_{i+1}$;
\item{(ii)} $\fork\bb {\ba_{n+1}}_{\ba_n}$.

Let $\ba_0=\lng \rng$. The induction step: as $\fork \bb {{M}} _
{\ba_n}$, by the choice of $\bb$, and the finite character of forking,
there is some finice $\bc \in M$ such that $\fork \bb \bc _
{\ba_n}$.
Let $\ba_{n+1}=\ba_n\concatenate \bc$. By monotonicity, $\fork \bb
{\ba_{n+1}}_{\ba_n}$. 

Now, we know that $\nonfork \bb {\cup_n\ba_n}_ {\cup_n\ba_n}$. By the
existence of non-forking extensions, we may assume that $\nonfork \bb
{{M}}_ {\cup_n\ba_n}$. 

We construct now by induction on $i<\l$ a continuous increasing
chain of models with the following properties:
\item{(1)} $N_i\sat T$ and $N_i\prec
N_{i+1}$. If $i$ is limit, 
$\ov N_i$ is the representation
$N_i=\bigcup_{j<i}N_j$. 
\item{(2)} For every $\eta\in {}^{<\om}i$,  strictly increasing,
$\ba_\eta\in N_i$. $\tp(\ba_{\eta\rest 0}\concatenate\ba_{\eta\rest
1}\concatenate\cdots\conc\ba_{\eta\rest{\lg\eta}})=
\tp(\ba_0\conc\ba_1\conc\cdots\conc\ba_{\lg\eta})$;
\item{(3)} If $\eta=\nu\conc \lng i\rng$, then $\nonforkin {\ba_\eta}
{N_i}_{\ba_\nu}^{N_{i+1}}$.
\item{(4)} If $i=\d\in S$ then there is some element $\bb\in 
N_{\d+1}$ such that $ \{\a_i\in c_\d:\forkin
{\bb}{N_{\a_{i+1}}}_{N_{\a_i}}^{N_{\d+1}}\}=Y$.

At the induction stage, when given $\nu\in {}^{<\om}i$ and increasing,
denoting by $\eta$ the sequence $\nu\conc i$, we should say who
$\ba_{\eta}$ is. There is an elementary mapping $h$ such
that for every $k<\lg\eta$, $h(\ba_k) = \ba_{\nu\rest k}$. Therefore
$h[\tp(\ba_{\lg\eta},\cup_{l<\lg\eta})]$ is a complete type over
$\bigcup_{l<\lg\eta} \ba_{\eta\rest l}$. By the existence of non
forking completion of a partial non forking type, there is some type
$p$ over $N_i$ which does not fork over $\bigcup
_{l<\lg\eta}\ba_{\eta \rest l}$. Let $\ba_\eta$ realize $p$ in
$N_{i+1}$.

In case $i=\d\in S^\l_0$ is as in (4), we should also take care of (4).

Let $Y(\d)=\lng \a^d_{i(n)}:i\in Y\rng$.
Let $\eta=\lng
\a^\d_{i(n)}: n\in \om \rng$, and let $h$ be an elementary mapping
such that $h(\ba_l)=\ba_{\eta\rest l}$. Then in $N_{\d+1}$ we 
add an element $\bb_\d$ which realizes $h(\tp(\bb,\cup_{l<\om}\ba_l))$
and $\nonforkin {\bb_\d} {N_\d} _ {\cup_{l<\om}\ba_{\eta\rest l}} ^
{N_{\d+1}}$ (due the existence of non forking extentions of types). We
have to show  

\ppro \lma Lemma: $\{\a_i:\forkin {\bb}{N_{\a_{i+1}}}_{N_{\a_i}}^{N_{\d+1}}\}=Y(\d)$.

We first need

\ppro \llma Lemma: If $\eta ^ \d(k)=\a_m$, then  $\nonfork
{\cup_{l<\om}\ba_{\eta\rest l}} {N_{\a_m}}_ {a_{\eta\rest
k}}$. 
\proof
By induction on $r,\quad k\le r< \om$, we see that $\nonfork
{\cup_{l\le r} \ba_{\eta\rest l}}  {N_{\a_m}} _
{\ba_{\eta}\rest k}$. 

$r=k$: $\nonfork{\cup_{l\le k}\ba_l}{N_{\a_m}}_ {\ba_{\eta\rest k}}$ is
trivial, as $\cup_{l\le k} \ba_{\eta\rest l}=\ba_{\eta\rest k}$. 

$r+1$: By the induction hypothesis, 
$$\nonfork {\cup _{l\le r} \ba_{\eta\rest l}} {N_{\a_m}} _
{\ba_{\eta^d\rest k}}\leqno {(a)}$$
\vskip1cm 

By the
construction,
$$\nonfork {\ba_{\eta\rest(r+1)}}
{N_{\eta(r)}}_{\ba_{\eta\rest r}}\leqno{ (b)}$$
\vskip1cm
Monotonicity gives 
$$\nonfork {\ba_{\eta\rest(r+1)}} {N_{\a_m}} _
{\ba_{\eta}\rest(r)}\leqno{\hfill (c)}$$ 
\vskip1cm

(a) and (c) give $$\nonfork {\ba_{\eta\rest(r+1)}}
{N_{\a_m}}_{\ba_{\eta\rest k}}$$

By the finite character of non forking Lemma \llma\ is proved.
\endproof{\llma}

Suppose now, first, that $\a_i\notin Y(\d)$. Let $\eta(k-
1)<\a_i<\eta(k)$.
We know that $\nonfork 
{\bb}{N_\d}_{\cup_{l<\om}\ba_{\eta}\rest l}$. So by monotonicity
$\nonfork {\bb}{N_{\a_{i+1}}}_{\cup_{l<\om}\ba_{\eta}\rest l}$. By
\llma, $\nonfork {\cup_{l<\om}\ba_{\eta}\rest
l}{N_{\eta(k)}}_{\ba_{\eta\rest(k)}}$. By monotonicity and the fact that
$\ba_{\eta\rest k}\in N_{\a_i}$ we get $\nonfork
{\cup_{l<\om}\ba_{\eta\rest l}}{N_{\a_{i+1}}}_{N_{\a_i}}$. By
transitivity of non forking, $\nonfork
{\bb_\d}{N_{\a_{i+1}}}_{N_{\a_i}}$, namely $\a_i\notin \inv_{\ov
N_{\d+1}}(\bb_\d,c_\d)$.

For the other direction: suppose that $\a_i\in Y(\d)$ and that
$\a_i=\eta(k)$.
We know that $\fork {\bb}{\ba_{\eta}\rest (k+1)}_{\ba_{\eta}\rest
k}$. Therefore by monotonicity, $\fork
{\bb}{N_{\a_{m+1}}}_{\ba_{{\eta}\rest
k}}$.
By \llma, as in the previous case, $\nonfork
{\bb}{N_{\a_{i}}}_{\ba_{\eta}\rest k}$. 
Suppose to the contrary that $\nonfork {\bb}{N_{\a_{i+1}}}_{N_{\a_i}}$.
Then by transitivity we get $\nonfork
{\bb}{N_{\a_{i+1}}}_{\ba_{\eta\rest k}}$ --- a
contradiction.\endproof{\lma,\construction}
 
We will need also

\ppro \secondConstruction Lemma: (the second construction Lemma)

Suppose $\l$ is uncountable regular, and $\ov C=\lng c_\b:\b<\l\rng$
is a weak club guessing sequence, such that for every $\b<\l$ the
order type of $c_\b$ is some fixed $\mu$ with $\cf\mu=\aleph_0$.
Suppose that $Y(*)\su \mu$ is given and of order type $\om_0$. There
is some model $M$ of $T$ with universe $\l$ and representation $\ov M$
such that $Y(*)\in P^*(\ov M, c_\b)$ for every $\b<\l$.

\proof The proof is essentially the same as this of \construction. The
only difference is in the construction: we add the witness not in
stage $\d+1$ but in stage $\b+1$, where $\sup
c_\b=\d$.\endproof{\secondConstruction}

\ppro \thirdConstruction Lemma:(the third construction lemma)
Suppose $T$ is a stable first order theory, $\cf\l=\l\ge |T|$, 
$\cf\k=\k<\k(T)$ and $\ov C, \ov P$ are as in \guessfacts\ (2).
Suppose $Y(*)\su \k$ is given. Then
there is a model $M\sat T$ of cardinality $\l$ and representation $\ov
M$ such that fore every $\d\in S$, $Y(*)\in P^*(\ov M,c_\d)$. 

\proof
We work in a monster model $\bbM$ and construct a sequence $\lng
\ba_\a:\a<\k\rng$ and an element $\bb$ such that $\ba_\a$ is
an infinite sequence, increasing with $\a$,  
namely $\ba_\a$ is
is a proper
initial segment of $\ba_\b$ whenever $\a<\b$ and $\fork
{\bb}{\ba_{\a+1}}_{\ba_\a}$ for all $\a<\k$. This is possible because
$\k<\k(T)$. Without loss of generality, $\nonfork
{\bb}{\cup_\a\ba_\a}_{\cup_\a\ba_\a}$. Let $Y(*)_\d\su c_\d$ for
$\d\in S$ be the isomorphic image of $Y(*)$ under the enumeration of
$c_\d$. We may assume, without loss of generality, that for every
$\a\in \nac c_\d$ for $\d\in S$, $Y(*)_\d\cap \a\in P_\a$. Construct
by induction on 
$\a<\l$ an elementary chain of models $M_\a$ with the following
properties:
\item{(1)}: for every $\eta\in P_\a$, $\eta\in [\a]^{<\k}$ there is a
sequence $\ov\ba_\eta$ such that $\ba_\eta(\b)\in N_\a$ and
$\tp(\cdots \conc\, \ba_\eta(\b)\,\conc\cdots)=\tp(\cdots
\conc\,\ba(\b)\,\conc\cdots)$. 
\item{(2)} If $\a=\d\in S$ then there is an element $\ba_\d\in
N_{\d+1}$ such that $\inv^*_{\ov M_{\d+1}}(\ba_\d,c_d)=Y(*)$. 

We let the reader verify that the analogs of \lma\ and \llma\ are
true.\endproof{\thirdConstruction}

\bigbreak
\neusection

{\bf \S\number\secno. The Main Results}
\ppro \countableTheorem Theorem:
Suppose $T$ is a complete, countable, stable but unsuperstable first
order theory, and that $\aleph_1<\l=\cf\l<2^{\aleph_0}$. Then
$\l\notin\univ(T)$. Furthermore, for every family $\{M_i\}_{i\in I}$,
$M_i\sat T$, $||M_i||=\l$ and $|I|<2^{\aleph_0}$ there is a model
$N\sat T$, $||N||=\l$ and $N$ is not elementarily embeddable into $M_i$
for all
$i\in I$. 

\proof

Clearly, it is enough to prove the ``furthermore'' part of the
theorem. Suppose  that $\{M_i\}_{i\in I} $ is a family of less then
$2^{\aleph_0} $ moels of $T$, each of cardinality $\l$. Let $\ov N_i$
represent $M_i$.  
Use \guessfacts\ 
part  1
to pick some club guessing sequence $\ov C$ on $S\su
\l$ with all $c_\d$ of order type $\om$. Pick some set $Y(*)\su \om$
such that $Y(*)\notin \bigcup_{i\in I,\d\in S}P^*(\ov N_i, c_\d)$. This
is possible, because the size of this union is smaller than
$2^{\aleph_0}$. Use the construction lemma to get a model $M$ of size
$\l$ and a  representation $\ov M$ 
such that for every $\d\in S$, $Y(*)\in P^* (\ov M,c_\d)$. Suppose
to the contrary that for some $i\in I$, $h:M\to M_i$ is an elementary
embedding.  By
\InvIsPreservedUnderEmbeddings, there is a club $E\su \l$ such that
for every $\d\in S$ such that $c_\d\su E$, $P(\ov M,c_\d)\su P(\ov
M_i,c_\d)$. Pick some $\d_0\in S_{E}$. So $Y(*)\in P^*(\ov
M,c_{\d_0})\su P^*(\ov M_i,c_{\d_0})$
-- a
contradiction to $Y(*)\notin \bigcup_{i\in I,\d\in S}P(\ov
M_i,c_\d)$.\endproof{\countableTheorem}

\ppro \uncountableTheorem Theorem:
Suppose $2^{\aleph_0}<\l=\cf\l<\l^{\aleph_0}$ and there are no $\mu_n$
such that $\l=(\Sigma \mu_n^{\aleph_0})^+$. Then if $T$ is a stable
unsuperstable theory, 
$|T|\le \l$, then $\l\notin\univ(T)$. Furthermore,  for every family
$\{M_i\}_{i\in I}$, 
$M_i\sat T$, $||M_i||=\l$ and $|I|<\l^{\aleph_0}$ there is a model
$M\sat T$, $||N||=\l$ such that $M$ is not elementarily embeddable
into $M_i$ for all 
$i\in I$. 

\proof
Again, the ``furthermore'' part is enough. 

Let $\mu$ be the least cardinal such that $\mu^{\aleph_0}> \l$. Since
$\l$ is uncountable and regular,
$\l^{\aleph_0}=\cup_{\a<\l}\a^{\aleph_0}$. If for every cardinal
$\k<\l$ $\k^{\aleph_0}=\k$, we should have had $\l^{\aleph_0}=\l$.
Therefore
$\mu$ is stricly smaller than $\l$. If $\cf\mu>{\aleph_0}$, then
$\mu^{\aleph_0}=\cup_{\a<\mu}\a^{\aleph_0}$. By the minimality of
$\mu$,
for every $\a<\mu$, $\a^{\aleph_0}\le \l$. This contradicts
$\mu^{\aleph_0}>\l$. We conclude that $\cf\mu={\aleph_0}$. Lastly, if
$\l=\mu^+$, then, $\mu$ being of cofinality $\om$, there would be
$\mu_n$ increasing to $\mu$ such that $\mu_n^{\aleph_0}<\mu$. This
contradicts the assumptions on $\l$.

Use \guessfacts\ 
part 3
to pick some weak club guessing sequence $\ov C=\lng
c_\b:\b<\l\rng$   with all $c_\d$ of order type $\mu$. Suppose to
the contrary that $\{M_i\}_{i\in I}$ is as stated above. 
By the assumption $\l<\mu^{\aleph_0}$, we can find some $Y(*)\su \mu$
of
order type $\om$ such that $Y(*)\notin \bigcup_{i\in I,\b<\l}P^*(\ov
M_i,c_\b)$. By  
the construction lemma there is some model $M$ and representation $\ov
M$
such that for every $\b<\l$, $Y(*)\in P^*(\ov M,c_\b)$.  Suppose to
the contrary that for some $i\in I$ there were an
elementary embedding $h:M\to M_i$. By \InvIsPreservedUnderEmbeddings\
there is a club $E\su \l$ such that if $c_\b\su E$ then $P^*(\ov
M,c_\b)\su P^*(\ov M_i,c_\b)$. As $\ov C$ is a weak club guessing
sequnce there is such a $c_\b$, and the contradiction to the choice of
$Y(*)$ follows as
before. \endproof{\uncountableTheorem}

\ppro \singTheorem Theorem:
Assume $T$ is first order complete countable stable  unsuperstable theory.
Suppose $\mu$ is singular, and there is some $\sigma<\mu$ and
$\k<\mu$ such that $\sigma^+<\k=\cf\k$ and
$\sigma^{\aleph_0}>\cov(\mu,\k^+,\k^+,\k)$, then there is no model of
$T$ in cardinality $\mu $ into which all models of $T$ of  cardinality
$\k$ are elementarily embeddable. In particular $\mu\notin\univ (T)$.

\proof We may assume that $\cf\sigma=\aleph_0$. 
Suppose to the contrary that $M\sat T$ is of cardinality $\mu$ and
that every $N\sat T$ of cardinality $\k$ is elementarily embeddable into it.
Without loss of generality the universe of $M$ is $\mu$. Let
$\theta\eqdef\cov(\mu,\k^+,\k^+,\k)$, and let $\lng B_i:i<\theta\rng$
demonstrate the definition of $\theta$. Without loss of generality,
each $B_i$ is the universe of some $M_i\prec M$ of cardinality $\k$.
 By \guessfacts\   part 3
pick some weak club guessing sequence $\ov C$
 with all 
$c_\b$ of order type $\sigma$. Pick a presentation $\ov M_i$ for every
$M_i$. Pick some $Y(*)\su \mu$ of order type $\om$ such that
$Y(*)\notin \bigcup_{i<\theta,\b<\k}P^*(\ov M_i,c_\d)$, and use
\construction\ to construct a model $N\sat T$ of cardinality $\k$ 
presetation $\ov N$
such that for every $\b<\k$, $Y(*)\in
P^*(\ov M,c_\b)$. For every $\b<\k$ there is some
element $\ba_\b$  such that $\inv^*_{\ov N}(\ba_\b,c_\b)= Y(*)$. 
Suppose that $h:N\to M$ is an elementary
embedding. There is some set of indices $X\su \theta$ such that
$|X|<\k$ and  $\ran
h\su \cup _{i\in X}B_i$.
Since
$\id^a(\ov C)$ is $\k$-complete, there is a set $S'\su S$, $S'\notin
\id^a(\ov C)$, and a fixed $i_0\in X$ such that $(\forall \d\in
S')(f(\ba_\d)\in B_{i_0})$.  Denote $B_{i_0}$ by $B$ for notational
simplicity, and let $L\prec M$ be the model with universe $B$. Use
\secondEmbedding\ to get the usual
contradiction.\endproof{\singTheorem}

\bigbreak      
\neusection

{\bf \S\number\secno.  Generalizations}

We wish now to generalise the discusstion of stable unsuperstable
theories --- namely those $T$ with $\k(T)=\aleph_1$ --- to stable
theories with $\k(T)$ arbitrary.

\ppro \general Theorem:
Suppose that $T$ is stable and that $\l\ge|T|$ is an uncountable
regular cardinal. Suppose that $\k<\k(T)$, and $\k^+<\l<2^\k$. Then
$\l\notin \univ(T)$. Furthermore, for every family $\{M_i\}_{i\in I}$
with $|I|<2^\k$ of models of $T$. each of cardinality $\l$, there is a
model $M\sat T$ of cardinality $\l$ which is not elementarily
embeddable into $M_i$ for all $i\in I$. 

\proof   
By \guessfacts\   part 2
there is a club guessing sequence $\ov C=\lng
c_\d:\d\in S\rng$ on some stationary set $S\su \l$ and a sequence $\ov
P=\lng P_\a:\a<\l\rng$ such that the order type of each $c_\d$ is
$\k$, for every $\a\in S$ and $\a\in \nac c_\d$, $c_\a\cap\a\in P_\a$,
and each $P_\a$ has cardinality $<\l$.
Pick a $Y(*)\su \k$ such that $Y(*)\notin \bigcup_{\d\in S,i\in
I}\inv^*(\ov M_i,c_\d)$ and use the third construction lemma to find a
model $M\sat T$ of cardinality $\l$ and a representation $\ov M$ such
that for every $\d\in S$, $Y(*)\in P^*(\ov M,c_\d)$. Suppose to the
contrary that there are $i\in I$ and an elementary embedding $h:M\to
M_i$. By \InvIsPreservedUnderEmbeddings\ and the fact that $\ov C$
guesses clubs we obtain the usual contradiction.\endproof{\general}

\ppro \shelah Theorem: 
Assume $\k = \cf(\k) < \k(T)$,    $\k \le \mu$,
$\mu^+ < \l = \cf( \l )  < \chi  <  \mu^\k$.
Suppose also that T is first order complete  
and $\k < \k(T)$.
Then there is no model $M$  of $T$ of cardianlity $\chi$
universal for 
models of $T$ of cardianlity $\l$.
\proof: similar.

\ppro \last Remark: This means that $(\l,1,\chi)  \notin \univ_t (T , \prec )$.
\bigskip
\neusection

{\bf\S\number\secno. A theory with a maximal universality spectrum}

In [KjSh 409], 5.5 it was shown that whenever $\l\in \univ(T)$, $T$ a
theory  having the 
strict order property, then there is a universal linear order in $\l$.
We prove now an analogous theorem for stable unsuperstable theories. 

\ppro \classDef Definition: for a cardinal $\k$, 
\item{(1)} $T_\k=Th(\lng {}^\k\om,E_\zeta\rng_{\zeta<\k})$ where
$\eta{E_\zeta}\nu\iff\eta\rest\zeta=\nu\rest\zeta$ 
( so $T$ is a first order complete theory of cardianlity 
$\k$  with $\k(T) = \k^+$,
and in fact is the canonical example of such a theory)
\item {(2)} $K_\k$ is the class of all trees of height $\k+1$. 
\item {(3)} $K^+_\k$ is the class of all trees of height $\k+1$  such
that   above every member thee is one of height $\k$

\ppro \equivFact Fact: $\univ (K_\k,<)\cap
(\k,\infty)=\univ(K^+_\k,<)\cap (\k,\infty)=
\univ(T_\k,\prec)\cap(\k,\infty)=\univ(T_\k,\le)$.
\proof Easy exercise.\endproof{\equivFact}

\ppro \minimality Theorem:
Suppose that $T$ is stable, $\k=\cf\k<\k(T)$, $\k\le\l$ and $|T|<\l\in \univ(T)$. 
Then $\l\in \univ(T_\k)$.
\ppro \afterlast Remark: Similarly for $\univ_p, \univ_t$

\proof  Without loss of generality, $|T|=\k$, for this may only
increase $\univ(T)$. So $|T|<\l$. 
Suppose that $N\sat T$ is universal in power $\l$. We define a model
$M$ which we shall prove to be universal in $\l$ for $K^+_\k$. By $\k<\k(T)$
we can find an 
element $\ba$ and an elementary chain $\lng M_i:i\le \l\rng$ such that
$\forks \ba {M_{i+1}} _ {M_i}$. 
 Let $M_\k^+$ be such that $M_\k\prec M_\k^+$ and
such that there is $\bI\su M_\k^+$, $|\bI|=\l$ and
$\bI$  an indiscrenible set based on $M_\k$,  i.e. $ \av ( \bI , \bbM)$
extends the type 
of $\ba$ over  $M_\k$
but   does not fork  
over $M_\k$.

The universe of $M$ will be $B=\{p\in 
S^1(N):p=\av(\bJ,N)$ for some $\bJ ,\;\bJ\su N,\;|\bJ|=
\l,\;\tp(\bJ)=\tp(\bI)\}$. 

\ppro \size Lemma: $|B|\le\l$
\proof Suppose to the contrary that there are $\l^+$ types $\lng
p_i:i<\l^+\rng$ and $\l^+$ indiscernible sets $\bJ_i\su N$, $|\bJ_i|=\l$ such that
$p_i=\av(\bJ_i,N)$. Pick a representation $\ov N=\lng N_\a:\a<\l\rng$
 of $N$ as an
elementary chain. For every $i<\l$ there is some $\a_i<\l$ such that
$|\bJ_i\cap N_{\a_i}|\ge \aleph_0$. Also, by \indisFacts\ (3) it
follows that  there is
some $c_i\in 
\bJ_i$ which realizes $\av(\bJ_i,N_{\a_i})$. By the pigeon hole
principle there are some
 $i<j<\l$ such that $\a_i=\a_j$ and
$c_i=c_j$. This contradics the fact that $p_i\not=p_j$ by \indisFacts\
(4). \endproof{\size}

By \indisFacts,(5), for every $p\in S(M_\k)$ and a finite set of formulas
$\Delta$ there 
is an infinite set of  indiscernibles 
$\bI\su M_\k$ such that $p=\av_\Delta(\bI,M_\k)$. By the stability of
$T$ and \indisFacts,(1),  there is some $n_\Delta$ such that for every
$\bJ\su\bI$ which satisfies $|\bJ|>2n_\Delta$,
$$(\forall \bb\in M_\k)(\forall\phi\in \Delta)(\phi(\ov x,\bb)\in p
\iff |\{\bc\in\bJ:\neg\phi(\bc,\bb)\}|\le n_\Delta\leqno{(*)}
$$

For every $\phi\in L$ there is a minimal $\a_\phi$ such that there is
a set $\bJ'_\phi\su M_{\a_\phi}$ of size $>2n_{\{\phi\}}$ which satisfies
$(*)$. Clearly, as $\bJ'_\phi$ is finite, $\a_\phi$ is a non limit
ordinal. By \indisFacts\ there is
 an infinite $\bJ_\phi\su M_{\a_\phi}$ with
$\av(\bJ,M_{\a_\phi})=p\rest M_{\a_\phi}$. By \indisFacts, $p\rest
M_{\a_\phi}=\av(\bJ,M_{\a_\phi})$. 

If $\sup\{\a_\phi:\phi\in L\}=\a^*<\k$, then $p$ were definable over
$M_{\a^*}$, and therefore, by \forkFacts, would not fork over
$M_{\a^*}$, contrary to its choice. We can, therefore, find a sequence
of formulas $\lng \phi_\zeta:\zeta<\k\rng$ with $\lng
\a_{\phi_\zeta}:\zeta<\k\rng$ strictly increasing. We shall assume, by
re-enumeration, if necessary, that $\a_{\phi_\zeta}=\zeta+1$. 

We define now the relations on our universe $B$. For every pair
$p_1,p_2\in B$ and $\phi_\zeta$ the following is an equivalence
relation: $p_1E^\zeta p_2\iff p_1\rest
\{\phi_\xi:\xi\le\zeta\}=p_2\rest\{\phi_\xi:\xi\le\zeta\}$.  Clearly,
these are $\k$ nested equivalence relations. We view  the structure we
defined as a tree of height $\k+1$ with no short branches,  
i.e.  is a member of$ K^+_\chi$.

To show that $M$ is universal, we will show that for every tree $S$ of
size $\l$ with all  branches of length $\k+1$ we
can find a model of $T$, $N_S$ of the same cardinality, such that the
elementary 
embedding of $N_S$ into the universal model $N$ will give an embedding
of $S$ into $M$. 

 We work by induction on $i\le \k$ and
for every $\eta\in {}^i\l$ construct an elementary embedding
$f_\eta:M_i\to \bbM$ with image $M_\eta$. We demand:
\item{(1)} $\nu\init \eta\imply f_\nu\su f_\eta$.
\item{(2)} for every $\eta\in {}^i\l$ and $\a<\l$, $\nonfork
{M_{\eta\conc\lng\a\rng}}{\bigcup\{M_{\nu\conc\lng\b\rng}:\nu\conc\lng\b\rng\in
{}^i\l,\;\nu\conc\lng\b\rng 
\not=\eta\conc\lng\a\rng\}\cup M_\eta}_{M_\eta}$

At limit $i$ we take unions. For $i+1$: $M_{\eta\conc\lng\a\rng}$
exists by \forkFacts.

For every $\eta\in {}^\k\l$ extend $f_\eta$ to $f_\eta^+:M_\e^+\to
\bbM$.

\ppro \ccc Claim: Suppose that $\nu\not=\eta\in {}^\k\l$ and that
$\zeta$ is the least such that $\eta(\zeta)\not= \nu(\zeta)$. Suppose
that $N\prec \bbM$ and that $M_\eta^+,M_\nu^+\su N$. Then
$\av(\bI_\eta,N)E^\xi\av(\bI_\nu,N)\iff \xi<\zeta$.

\proof Let $\xi<\k$. Let $f\in\aut(\bbM)$ map $M_\eta^+$ onto
$M_\nu^+$. For simplicity we assume that $f_\eta\rest M_\zeta=f_\nu\rest
M_\zeta=id$. We know that $\nonforkin
{M_\eta}{M_\nu}_{M_\zeta}^{\bbM}$. 

First case: $\zeta>\xi$.  $(*)$ gives a definition of $\av_{\phi_\xi}(\bI,N)$
with set of parameters $\bJ'_{\phi_\xi}$. In $\bbM$ this definition
gives, with respective sets of parameters $\bI_\eta,\bI_\nu$, the
types $p_1,p_2$, which extend, respectively, $\av_{\phi_\xi}(\bI_\eta,M_\eta^+)$
and $\av_{\phi_\xi}(\bI_\nu,M_\nu^+)$. Let $\bI_\eta$ be $\bI_1$ and
let $\bI_\nu$ be $\bI_2$.

\ppro \crux Fact: for $l=1,2$,
$\av_{\phi_\xi}(\bI_l,\bbM)=\av_{\phi_\xi}(\bJ_{\phi_\xi},\bbM)$.

\proof Suppose to the contary that $\bc\in \bbM$ demonstraes
otherwise. Then by \indisFacts, there is someh is 
 $\bI'_l\su \bI_l$
of size $<\k(T)$ such that the set $(\bI_l\sm \bI'_l)\cup \{\bc\}$ is
independent over $M_l$.  Thereforee $\nonfork \bc {M_l + \cup(\bI_l\sm
\bI'_l)}_{M_l}$. By \forkFacts (4), the type of $\bc$ is finitely
satisfiable over $M_l$. There is finite information saying that
$\phi(-,\bc)$ behaves ddifferently in
$\av_{\phi_\xi}(\bJ_{\phi_\xi},\bbM)$ than in
$\av|{\phi_\xi}(\bI_l,\bbM)$. So there is a counterexaample inside
$M_l$ --- a contradiction. \endproof{crux}

By this fact we conclude that
$\av(\bI_\eta,\bbM)E^\xi\av(\bI_\nu,\bbM)$.

Second case: $\xi\ge \zeta$. We extend $\av(\bI_\eta,M_\eta)$ to a
non-forking extention $p\in S(\bbM)$. So $\nonfork p
{M_\nu}_{M_\zeta}$. In particular $\nonfork {P\rest
\phi_\xi}{M_\nu}_{M_\zeta}$. Therefore there is some $\bJ_\zeta'$ as
in $(*)$ --- contradiction to $\a_\xi=\xi+1$. \endproof{ccc}

Suppose now that $S$ is a given tree in $K^+_\k$ of size $\l$. Without
losssof generality, $S<{}^\k\l$. Pick a model $N_S\prec \bbM$ such that
for every $\eta\in S$, $M_\eta\su N$ and such that $||N||\le \l$. An
elementary embedding of $N_S$ into $N$ easily gives an embedding of
$S$ into $M$. \endproof{\minimality}

\bigbreak

\centerline{\bf REFERENCES}

\item{[KjSh 409]} M. Kojman and S. Shelah, 
Non-existence of universal orders in many cardinals, accepted to the
{\bf
Journal of Symbolic Logic}

\item{[Sh-c]} Sahron Shelah,  {\bf Classification theory:  and the
number of  non- 
isomorphic models},
revised, {\it North Holland Publ. Co.}, Studies in Logic and the 
Foundation of 
Math vol. 92, 1990, 705 + xxxiv.

\item{[Sh-g]} Saharon Shelah, {\bf Cardinal Arithemetic}, to appeare.

\item{[Sh 100]} Saharon Shelah, Independence results, {\bf Journal of
Symbolic Logic} {\bf 45} No. 3 (1980) pp. 563--573.

\item{[Sh 175a]} Saharon Shelah, Universal graphs without CH:
revisited, {\bf Israel Journal of Math} {\bf 70} (1990) pp. 69--81

\item{[Sh 420]} Saharon Shelah, Cardinal arithmetic, to appear in 
Proceedings  of the Banff conference,
April 1991.

\bye